\documentclass{amsart}
\usepackage{a4}
\usepackage{amssymb}
\def\FF{{\mathbb F}}
\def\squareforqed{\hbox{\rlap{$\sqcap$}$\sqcup$}}
\def\qed{\ifmmode\squareforqed\else{\unskip\nobreak\hfil
\penalty50\hskip1em\null\nobreak\hfil\squareforqed
\parfillskip=0pt\finalhyphendemerits=0\endgraf}\fi\medskip}

\newcommand{\udot}{{}^{\textstyle .}}

\newcommand{\PSL}{\mathrm{PSL}}

\newcommand{\Co}{\mathrm{Co}}

\newcommand{\Mat}{\mathrm M}
\newcommand{\Suz}{\mathrm{Suz}}

\title{The Monster and black-box groups}
\author{Robert A. Wilson} 
\address{School of Mathematical Sciences, Queen Mary University of London,
Mile End Road, London E1 4NS, U.K.}
\date{First partial draft 30/07/2013; this version 17/10/2013}

\begin{document}
\maketitle

\centerline{\em Dedicated to the memory of \'Akos Seress}

\begin{abstract}
We discuss ways in which the black-box model for computation
is or is not applicable to the Monster sporadic simple group. Conversely,
we consider whether methods of computation in the Monster
can be generalised to other situations, for example to
groups of `cross-characteristic' type.
\end{abstract}

\section{Introduction}
\label{intro}
The concept of black-box group was introduced by Babai and
Szemer\'edi \cite{BSz} as an abstraction and generalisation of the ideas
of `matrix group' and `permutation group', for the purposes of algorithm
design and complexity analysis. In a black-box group, the elements are represented
(not necessarily uniquely) by bit-strings of some fixed length $n$,
and there are `black boxes' that perform the three operations
of group multiplication, inversion, and testing whether a given bit-string represents
the identity element, each in a specified
maximum amount of time.

The Monster is the largest of the $26$ sporadic simple groups, and the only one
for which no matrix or permutation representation is small enough for
naive computation to be effective (yet).
For practical computations we generally use 
instead the computer construction described in \cite{2loccon},
in which the Monster is generated by a subgroup $G=C(z)=\langle a,b\rangle \cong
2^{1+24}\udot \Co_1$, together with a `triality element' $T$ of order $3$,
which centralizes
a subgroup $2^{11}\udot \Mat_{24}$ of $\langle a,b\rangle$.

Thus an element of the Monster is represented as a word in $T$ and elements
of $G$. Multiplication is concatenation of words, combined with reduction by
$T^3=1$ and multiplying together any contiguous elements of $G$.
Inversion can be implemented as reversal of a word,
followed by 
replacing each element of $G$ by its inverse, and each occurrence of $T$ by $T^{-1}$.
And there is a quick and straightforward test for whether
a given word represents the identity element.

It is a natural question to ask, to what extent this form of computation is covered by
the black-box paradigm, or whether a different model is required. 
At first glance, the Monster appears  to conform to the black-box paradigm,
except that 
\begin{itemize}
\item we have not been given an effective bound on the number of bits
required to represent any element of the group;
\item the time taken for both the inverse operation and the identity test is 
proportional to the length of the word, so is potentially unbounded.
\end{itemize}

So, various questions arise as to the extent to which the insights gained from
the black-box approach are applicable to the very real problem of
actually computing anything interesting in the Monster. Issues of
complexity, of course, do not arise, but issues of efficiency are of
paramount importance. Moreover, since the Monster is so large,
efficiency questions do almost look like complexity questions.

In this paper we attempt to analyse computational questions about the Monster
from this point of view. Conversely, we ask to what extent insights gained
from computing in the Monster can be applied more generally, to problems
which are usually considered in the black-box model.

\section{Can we bound the length of the words?}
\label{wordshort}
Unless we have effective methods of shortening words,
the Monster fails the black-box principle that there should be an effective upper bound on
the number of bits needed to represent an element of the group.
For these purposes, a word is taken to be an alternating sequence 
of elements $g_i\in G$, and $T_i=T^{\pm 1}$, and we define the \emph{length} of a word
to be the total number of $g_i$ plus the total number of $T_i$.
Note that this is a different definition from \cite{L229}, where only the
number of $T_i$ is counted.
Multiplication of elements is, in the worst case, simply concatenation
of words, and therefore normal paradigms of computation
lead to exponential growth in the length of the words, and consequently
exponential growth in the number of bits used in the representation
of elements. 

Some extremely useful
methods for shortening certain words in certain circumstances
are described in \cite{L229,L259,L241} and other papers cited below.
Combining these in ways suggested by Ryba's 
constructive 
membership testing algorithm \cite{Ryba} leads to an effective
method for shortening any word of length greater than some
reasonable bound (which
can be taken to be as small as $17$ if we wish).
However, this is a randomised method, which has a small but non-zero
probability of failure. Moreover, from a theoretical point of view,
it is not clear whether repeated attempts are sufficiently independent
to ensure eventual success in practice, or whether there actually exist 
elements
of the Monster for which the method is guaranteed to fail.
(This latter possibility, however, seems very unlikely, and even if it happens,
small changes to the method should eliminate it, at the cost of an increase
in the bound of $17$.)

The first basic  building block 
was introduced in \cite{L229}, and
is a method of taking a word which represents an
element of $C(z)$, and writing it in canonical form as a word of length $1$.
We may then make a recursive application of Ryba's algorithm
(or some other constructive membership test) in $2^{1+24}\udot\Co_1$
to find a word in $a$ and $b$ which is equal to the given Monster element. In fact,
Ryba's algorithm does not work very well in groups with large normal $2$-subgroups,
so we make some modifications, as described below.

The second basic building block of the method
was introduced in \cite{L259}, where it  was called `changing post'.
What this means is, given  any $2B$-element $t$ which centralizes $z$, 
finding a word (of length at most $4$) which
conjugates $t$ into $z$.
This in itself does not shorten any words: indeed, it lengthens them.
The original word of length $1$ for $t$ is turned
into a word of length (typically) $9$.

The technical details of these two processes are described in the next section.
Now we explain how they are combined into a word-shortening algorithm.
 
Given any word $W$ representing an element of the Monster,
\begin{enumerate}
\item
take (random) elements $g\in G$ until the word $W':=Wg$
is an element 
which powers up to an involution $t$ in
class $2B$;
\item conjugate $t$ by (random) $c_0\in G$ until $t^{c_0}z$ has even order and
powers up to a $2B$-element $y$;
\item find the word of length $1$ for $y$;
\item find a word $c_1T^{\gamma_1}c_2T^{\gamma_2}$ 
of length $4$ which conjugates $y$ to $z$;
\item find the word of length $1$ for $t^c$, 
where $c=c_0c_1T^{\gamma_1}c_2T^{\gamma_2}$;
\item find a word $d=d_1T^{\delta_1}d_2T^{\delta_2}$ 
of length $4$ which conjugates $t^{c}$ to $z$;
\item find the word of length $1$ for $W'':=W'^{cd}$.
\end{enumerate}
Finally we have
\begin{eqnarray*}
W&=&W'g^{-1}\cr
&=&cdW''d^{-1}c^{-1}g^{-1}\cr
&=& c_0c_1T^{\gamma_1}c_2T^{\gamma_2}d_1T^{\delta_1}d_2T^{\delta_2}
W''T^{-\delta_2}d_2^{-1}T^{-\delta_1}d_1^{-1}
T^{-\gamma_2}c_2^{-1}T^{-\gamma_1}c_1^{-1}c_0^{-1}g^{-1}
\end{eqnarray*}
Since we can compute products within $G$, the total length of the word for $W$
is $17$. In certain (unlikely) circumstances, the words for $c$ and/or $d$ may be
shorter, in which case the word for $W$ is correspondingly shorter than $17$.

\section{The technical details}
\label{tech}
This section can be skipped by those readers who only want an overview,
or an understanding of the general principles.

\subsection{The underlying module}
The generators, $T$ and elements of $G$, for the Monster, are stored in a format
whereby their actions on the underlying module, of dimension $196883$ 
over $\mathbb F_3$, can be computed. This is the foundation for both the
identity test and the order oracle. Two vectors have been pre-computed,
whose joint stabilizer is proved to be trivial. Hence a word represents the
identity element if and only if it fixes these two vectors.
(If one is prepared to make do with a Monte Carlo algorithm, one can
halve the cost by taking one random vector instead.)

\subsection{Changing post}
There are just five classes of involutions in $C(z)$ which lie in Monster class $2B$.
For a representative $x$ of each class we pre-compute a word which conjugates
$x_i$ to $z$, as follows.
\begin{enumerate}
\item In the case $x_1=z$, there is nothing to do.   
\item In the case $x_2\in 2^{1+24}$, we may choose $x_2$ such that $x_2^T=z$.
\item Otherwise, $x_i$ maps to an involution in the quotient $\Co_1$, and
we may choose $x_i$ such that $x_i^T\in 2^{1+24}$. In each of the
three cases we search for an element $y_i$ of $G$ which conjugates
$x_i^T$ to the canonical representative $x_2$ of this class (see below for details). 
\end{enumerate}

If now $x$ is any $2B$-element in $C(z)$, conjugate to $x_i$, say, we search for
an element which does this conjugation.
For $i=3,4,5$, we perform the conjugation in $\Co_1$ first, and then
conjugate by suitable elements of $2^{1+24}$ as necessary afterwards.
Since any involution centralizes at least $2^{12}$ in $2^{1+24}$, even
an exhaustive search is not impossible.
For $i=2$, we adopt a randomised approach, and fingerprint around
$1000$ conjugates of each of $x$ and $x_2$. Sorting and merging
the two lists
of fingerprints we easily find a match, and read off an element of $G$ which
conjugates $x$ to $x_2$.

In all cases we now have a word $c_1Ty_iT$, or $c_1T$, or the empty word,
which conjugates $x$ to $z$. In fact for technical reasons it is easier to allow
also the possibility of using $T^{-1}$ rather than $T$.

\subsection{Computing words of length $1$ for elements which centralize $z$}
We first work in the quotient $\Co_1$ of $C(z)$, so that Ryba's algorithm
(or some other constructive membership test) can be applied directly.
It is straightforward, if somewhat technical, to obtain elements of $2\udot \Co_1$
as $24\times 24$ matrices over $\FF_3$, corresponding in pairs (modulo sign)
to elements of the quotient of $2^{1+24}\udot \Co_1$ by the normal $2$-subgroup.

This process can be carried out for
any element of the Monster which commutes with the central
involution $z$ of $2^{1+24}\Co_1$, even if it is only given as a word in the 
generators of the Monster. We just have to compute the images of a carefully
selected set of $24$ coordinate vectors, and extract another 
(not necessarily the same) carefully
selected set of $24$ coordinates from the answer.

Now we need to lift to $2^{1+24}\Co_1$.
Suppose that $w$ is the (long) word, and $x$ is the element
given as a word in $a$ and $b$. Then we know that $wx^{-1}$ is an element
of the group $2^{1+24}$. It so happens that there is an easy
constructive membership test in $2^{1+24}$. Thus we have words in
$a$ and $b$ for both $x$ and $wx^{-1}$, and we combine them to get a
word in $a$ and $b$ for $w$.
%

\section{Is the Monster a black-box group?}
\label{M=BB}
Of course, this is a meaningless question. Or at best, it is a philosophical
question, not a mathematical one. (A black-box group, after all, is one which
in principle one knows nothing about.) One can perhaps best interpret the question at a
purely phenomenological level, and ask whether black-box algorithms
(a) work at all, or (b) are effective, in the given computational
environment for the Monster.

With the word-shortening method described in the previous sections,
we have some bounds on the number of bits required to represent
any element of the Monster, and the time required for each operation.
In the current implementation, $n$ is around $6\times10^{8}$.
(This compares with about $7\times 10^{10}$ for the underlying matrices,
or $4\times 10^{10}$ for the minimal representation.)
The time required for the identity test is around $5$ seconds, and that
for inversion is about $1$ minute.

The multiplication algorithm envisaged here has not been implemented,
but, assuming that typical elements will require close to the maximum
word length of $17$, all the work is in shortening the resulting word of
length around $33$. Currently every shortened word we require
is made by hand,
and takes a day or two to make. If the method were automated and 
efficiently implemented I guess it
would take an hour or two.
(For comparison, the time taken to multiply two $196882\times196882$
matrices over $\mathbb F_3$ on the same system
would be more than a week.)

So, on the face of it, this would make the Monster into a black-box group
according to the official definition. 
(Actually,
this is not quite true, because the black box for element multiplication is
now only a Las Vegas algorithm: it may report failure instead.)
Black-box algorithms can be used,
although the very high cost of multiplication is a barrier.
Indeed, the computations that
are required for practical problems like determining the maximal subgroups,
require hundreds of thousands of multiplications, and therefore would
take years, as opposed to the small number of days taken by the
calculations we have actually done. 

For these reasons, then, we conclude that the
black-box model is still not a very useful model for the Monster.

\section{Compare and constrast}
\label{compare}
It is clear by now that,
at least in some respects, we are in a better position than
in the black-box situation. 
This must be so, for black-box algorithms alone could not
achieve the results that are described in the papers we cite below.
To compensate for the lack of a generic multiplication algorithm,
various other techniques are available.

Most importantly, a fantastically efficient order oracle is available.
Since every element has order at most $119$, computing the order of a word 
can be done in at most $119$ times the time taken to test whether it is
the identity element, so around $10$ minutes for a word of length $17$.

Philosophically, the black-box model is a socialist model: all elements
of the group are treated equally. 
The basic operation is multiplication of (arbitrary) group elements,
and complexity is measured in terms of the (worst case) time taken for this
basic operation.
An order oracle is often assumed, and is generally taken to be
at least as expensive as multiplication of elements.

But the model of computation in
the Monster is an elitist one: the elements of the subgroup $G$ are highly
favoured, because they can be multiplied together in about
$5$ seconds with no increase in word length. Multiplication of
\emph{arbitrary} group elements, on the other hand, is too
expensive for indiscriminate use. Practical computations
tend to be dominated either by computations in $G$, or by the order oracle,
depending on the context.

We should also consider to what extent matrix invariants are available
in the Monster. Each word which represents an element of the Monster
can in principle be converted into a $196882\times 196882$ matrix over $\mathbb F_3$.
The time taken for this operation is of the order of a few hours per letter of
the word. Thus a trace can be computed in around a couple of days.
This is unlikely ever to be cost-effective, however.
Other matrix operations are in general not available
in an effective manner. We know of no efficient way of calculating
 the characteristic polynomial, or  the Jordan block structure,
for example.

\section{The Monster as an infinite group}
\label{Minf}
 The ideas in this
section were expressed to me by Sasha Borovik.
It has often been remarked that the Monster 
is so large that it is `morally' an infinite group, or even
`bigger' than many infinite groups.
As he put it to me, the Monster is \emph{de jure} finite
but \emph{de facto} infinite. The same may be said of black-box groups,
if in a somewhat different way.

The difference, as Borovik explained it, is that black-box groups have invariant
probability measures, whereas the Monster does not.
Therefore black-box groups behave like compact Lie groups,
while
the Monster 
behaves 
more like
an HNN-extension or free amalgamated product.

Indeed, the form of the Monster construction we are using even
looks like an HNN-extension.
It is not actually \emph{free}, of course, in the sense that an HNN-extension
is free, but for the purposes
of many computations, it might as well be.
Computing in the Monster is very like computing in an HNN-extension, in that
very little can be done except in conjugates of the base subgroup.

\section{The Pacific Island model}
\label{Pacific}
To use a geographical analogy, the group $G$ is a small island where
productive work can be done, in a Monstrously vast ocean of other
group elements most of which are apparently useless. Occasionally
one has to fish in this vast ocean for elements outside $G$ that
will perform useful functions. (Perhaps $T$ is the boat that enables us
to travel on these fishing trips?---$T$ stands for `travel' as well as `triality'.) 
And when even that fails, one has
to navigate to distant islands and trade for the elements one requires
($T$ also stands for `trade').
Typically, this may involve scouring an entire unfamiliar island
is search of the elusive prize. (Actually, the Pacific Ocean is too small, and the
islands too large, for this analogy. The number of elements in the Monster is about
100 million times the number of water molecules in the Pacific Ocean, whereas the
number of elements in $G$, if converted into silica molecules, will give you 
not an island, but a bucket of sand.)

What are the essential ingredients of this model of computation?
We seem to need the following:
\begin{itemize}
\item black boxes to perform group operations in $G=C(z)$;
\item an order oracle for (short) words;
\item an oracle to solve constructively the conjugacy problem in the
conjugacy class of $z$;
\item an oracle for constructive membership testing in $G$.
\end{itemize}
These are listed in approximate order of cost in the Monster,
from cheapest to most expensive. 

The first and last are not
really computations in the whole group, but only in the subgroup $G$,
so should perhaps be taken as read. The third we have already shown
how to do, given the other three. This reduces the essential
requirements to one thing only: 
namely, an order
oracle, which is much cheaper than multiplication of elements.

To summarize, in the Pacific Island model of group computation, we are given
\begin{itemize}
\item an involution centralizer $G=C(z)$, in which all problems can be solved,
\item one extra generator $T$ (which probably should satisify $zz^T=1$),
and 
\item an order oracle for words of bounded length, 
\end{itemize}
and nothing else.
%

\section{Pacific Island, or Pacific Ocean?}
\label{ocean}
As currently practiced, computation in the Monster is largely
restricted to working on one particular island, that is the subgroup $G$,
with the occasional
fishing or trading trip to collect additional elements to perform particular
functions.
This is not the only possible way to work, however. Often we want to
work on other islands, i.e. compute in different subgroups. 

Occasionally this has been done in the Monster. A collection of words
is found for elements generating a (usually maximal, or close to
maximal) subgroup. The representation of this subgroup on the
underlying $196882$-dimensional module is investigated via a
specially designed form of `condensation', in order to construct
explicitly the action on a suitably small invariant subspace.
Then this subgroup can be explored in the usual way as a matrix group.

\section{Is the Pacific Island model useful in other contexts?}
\label{PImodel}
It is of course well-known that the black-box model, while extremely useful,
does not capture every important aspect of computation in finite groups.
It is commonplace to use other information if it is available,
for example traces of matrices, numbers of fixed points of permutations,
order oracles, and so on. 
But most of this extra information is still used in a socialist paradigm:
all group elements are treated equally.

In practice, however, many modern algorithms, for matrix groups
in particular, recurse to a subgroup,
often an involution centralizer, as quickly as possible. 
Efficient implementations will generally convert elements of the subgroup
into a form where computations proceed much more quickly.
In this scenario, the black-box model seems less applicable than the
Pacific Island model alluded to above.

The black-box approach is most useful in the beginning stages of
an investigation, when we have essentially no knowledge about
the group under discussion. But in the later stages, we generally
have a lot of knowledge, and often even know the isomorphism
type of the group (up to a certain probability of error). This again
favours the Pacific Island model.

However, there remains the important question as to which
of the additional operations that are practical in the Monster
remain practical and efficient in these other contexts. 
The crucial issue is whether or not there is a fast order oracle.

On the face of it, it seems hard to imagine many contexts
in which an order oracle is much faster than 
a single multiplication! What makes it work in the Monster is the fact that
the degree of the representation is much bigger than the largest element order.
This is a phenomenon associated with sporadic groups, or
cross-characteristic representations.

It is probably not a coincidence, that these are situations in which the
black-box model does not have a great deal to tell us. To make a sweeping
generalisation, in this situation the input data is so large compared to
the order of the group, that almost
all algorithms are more-or-less linear. But that does not necessarily mean that
all problems are soluble \emph{in practice}: just as in the Monster.
 So maybe 
the Pacific Island model can tell us something useful
about how to perform such calculations?

\section{Pacific island algorithms}
\label{PIalg}
The word-shortening algorithm described in Section~\ref{wordshort} is in effect
a Monte Carlo `constructive membership testing' algorithm. It takes an element
of the Monster, and writes it as a word in the standard generators. If the input
is not in fact an element of the Monster, then the algorithm probably fails,
but may just give the wrong answer.

Another problem one might wish to solve is the constructive recognition problem:
given another group which is claimed to be the Monster, can we produce an
isomorphism with the standard copy? This problem comes in various flavours,
depending on how much information is compatible between the two copies.

Suppose first that we are given copies of $G$ in both groups. Then, by assumption,
we can find isomorphisms between these copies of $G$. If $T$ is also given in both groups,
then we want to adjust these isomorphisms so that they map $T$ to the right place.
In each copy of $G$ we can find the centralizer of the relevant copy of $T$, and then
inside $G$ we can conjugate one of these centralizers to the other. It remains only
to identify which is which of the 8 elements of order $3$ in $\langle z,T\rangle\cong A_4$,
which can be achieved by computing orders of a few random elements.

A more likely scenario, which has actually occurred in practice, is that
$G$ is not immediately available in one of the copies of the Monster.
For example, we might have the `mod $2$' construction, in which, instead of $G$, a
subgroup $H\cong 3^{1+12}\udot 2\udot\Suz{:}2$ is the island in which
one can work. In \cite{ML227} we have given a method of obtaining
generators for $H$ from those of $G$, and then constructive recognition
within $H$ allows us to find the standard generators. Moreover, the extra generator
for the $H$-type Monster can be found within $G$, although this has
not actually been done yet.

The isomorphism in the other direction may be obtained by a similar process.
First, there is an involution whose centralizer in $H$ is
$6\udot\Suz{:}2$. The full centralizer of this involution in the Monster is
$2^{1+24}\udot\Co_1$, and involution centralizers can be found in the
Monster either by Bray's algorithm, or by some more subtle technique.
On the other hand, $T$ does not lie in $H$, so cannot be found quite
so easily.

In a more general situation, one might imagine that $G=C(z)$ and $H=C(y)$
are both involution centralizers, in some large `cross-characteristic'
group. In this case we would probably want to choose
the isomorphisms so that $z$ and $y$ commute with each other. If it is possible
also to choose $T$ commuting with $y$, then the above method will
construct a suitable isomorphism.

On the other hand, these are not really the type of calculations which are
most often going to arise in the Pacific Island situation. More likely, we already
know what the group is, and we are interested in computing particular subgroups,
such as centralizers and normalizers, Sylow subgroups, and the like.
The references below give many examples of calculations of
this type in the Monster.
In this context, it would be normal to assume that we are working in a `standard copy'
of the group. This allows us the luxury of pre-computing a great deal of the
structure. But it also rules out the cross-characteristic groups, and leaves us only with
sporadic groups. 

\section{Conclusion}
By comparing and constrasting (a) the black-box model of computation in finite groups,
and (b) the Pacific Island model of practical computation in the Monster, we have seen
that, in fact, there is surprisingly little overlap between the two, and surprisingly
little application of either method in the other context.

There have been one or two major influences from the sporadic context to the
black-box context: especially, the emphasis on computing involution centralizers;
and arising from that, Ryba's constructive membership algorithm. 
Moreover, black-box methods often involve a reduction to simple groups,
in which case the sporadic groups have to be dealt with somehow.
In the other
direction, subgroups of sporadic groups can be efficiently investigated using
general-purpose black-box algorithms. 

But beyond these influences, the two
inhabit rather different worlds. 
They complement each other, and both are necessary for a fully functional
toolkit for computation in finite groups.

\end{document}